\documentclass[12pt,a4paper]{article}

\usepackage{amsmath,amssymb,amsfonts}
\usepackage{epsfig}

\usepackage[latin1]{inputenc}
\usepackage[T1]{fontenc}
\usepackage{ae}
\usepackage[english]{babel}

\usepackage{fancyhdr}
\usepackage{lastpage}
\usepackage{geometry}
\numberwithin{equation}{section}

\newcounter{item}

\renewcommand{\theitem}{\arabic{section}.\arabic{item}}

\newcommand{\cc}{\setcounter{equation}{0}}

\newenvironment{theo}[1]{
\setcounter{item}{\value{equation}}
\addtocounter{equation}{1}
\refstepcounter{item}
\par\addvspace{\bigskipamount}
\indent {\bf \theitem.\hspace{1em}Theorem#1.} \sl }
{\par\addvspace{\bigskipamount}
}

\newenvironment{pf}{
\par\addvspace{-\smallskipamount}
\indent {\bf Proof.}$\,\ $ } { $\Box$
\par\addvspace{\bigskipamount}
}

\newenvironment{lem}{
\setcounter{item}{\value{equation}} \addtocounter{equation}{1}
\refstepcounter{item}
\par\addvspace{\bigskipamount}
\indent {\bf \theitem.\hspace{1em}Lemma.\,\ }  \sl }
{\par\addvspace{\bigskipamount}
}

\newenvironment{cor}{
\setcounter{item}{\value{equation}} \addtocounter{equation}{1}
\refstepcounter{item}
\par\addvspace{\bigskipamount}
\indent {\bf \theitem.\hspace{1em}Corollary.\,\ }  \sl }
{\par\addvspace{\bigskipamount}
}

\newenvironment{nonsec}{\bf
\setcounter{item}{\value{equation}}\addtocounter{equation}{1}
\refstepcounter{item}
\par\addvspace{\bigskipamount}
\indent \theitem.\hspace{1em}\ignorespaces }
{\unskip .\ \ \ }

\newcommand{\comment}[1]{}
\newcounter{minutes}
\divide\time by 60
\newcounter{hours}
\multiply\time by 60
\addtocounter{minutes}{-\time}

\begin{document}

\begin{center}
{\Large \bf On harmonic quasiconformal  quasi-isometries}
\end{center}
\medskip

\begin{center}
{\large \bf  M. Mateljevi\'c and M.~Vuorinen}
\end{center}

\bigskip
\begin{center}
{\tiny \texttt{FILE:~\jobname .tex,
        printed on: \number\year-\number\month-\number\day,
         \thehours.\ifnum\theminutes<10{0}\fi\theminutes}}

\end{center}

\medskip

{{\bf Abstract.} The purpose of this paper is to explore conditions
which guarantee Lipschitz-continuity of harmonic maps w.r.t.
quasihyperbolic metrics. For instance,  we prove that harmonic
quasiconformal maps are Lipschitz w.r.t. quasihyperbolic metrics.}

{{\bf 2000 Mathematics Subject Classification.} Primary 30C85.
Secondary 30C65.}

\bigskip

\section{Introduction}{}

\smallskip

Let $G \subset {\mathbb R}^2$ be a domain and let $f: G \to {\mathbb
R}^2, f= (f_1, f_2),$ be a harmonic mapping. This means that $f$ is
a map from $G$ into $ {\mathbb R}^2$  and both $f_1 $ and $f_2$ are harmonic
functions, i.e. solutions of the two-dimensional Laplace equation
\begin{equation}
\label{eq1.1}
\Delta  u = 0 \, .
\end{equation}
The Cauchy-Riemann equations, which characterize analytic functions,
no longer hold for harmonic mappings and therefore these mappings
are not analytic.
Intensive studies during the past two decades show
that much of the classical function theory
can be generalized to harmonic mappings
(see the recent book of Duren \cite{D} and the survey of Bshouty
and Hengartner \cite{BH}). The purpose of this paper is to
continue the study of the subclass of quasiconformal and harmonic
mappings, introduced by Martio in \cite{ma} and further studied
for example in
\cite{ma3,ma4,ma.koebe.ff,mk,ps,ka,akm,mak,kam2,kam3, kp, dkan}. The above
definition of a harmonic mapping extends in a natural way to the
case of vector-valued mappings $f: G \to {\mathbb R}^n, f=
(f_1,\dots, f_n),$ defined on a domain $ G \subset {\mathbb R}^n,
n \ge 2.$

We first recall the classical Schwarz lemma  for the unit disk
${\mathbb D} =\{ z \in {\mathbb C} : |z|<1  \}\,$:

\begin{lem}{} \label{ClassSLemma}
  Let $f: {\mathbb D} \to {\mathbb D} $ be an analytic function
with $f(0)=0\,.$ Then $|f(z)| \le |z|$, $z \in {\mathbb D} \,.$
\end{lem}

For the case of harmonic mappings this lemma has the following
counterpart.




\begin{lem}{{\rm (\cite{he}, \cite[p. 77]{D})}} \label{SLemma}
Let $f: {\mathbb D} \to {\mathbb D} $ be a harmonic mapping
with $f(0)=0\,.$ Then $|f(z)| \le \frac{4}{\pi} {\mathrm{
tan}}^{-1}|z|$ and this inequality is sharp for each point $z \in
{\mathbb D} \,.$
\end{lem}

The classical Schwarz lemma is one of the cornerstones of geometric
function theory and it also has a counterpart for quasiconformal
maps (\cite{A,lv,qiuvavu,vu2}). Both for analytic functions
and for quasiconformal mappings it has a form that is conformally
invariant under conformal automorphisms of $  {\mathbb D} \,.$

In the case of harmonic mappings this invariance is no longer
true. In general, if $\varphi : {\mathbb D} \to  {\mathbb D} $ is
a conformal automorphism and $f :  {\mathbb D} \to  {\mathbb D}  $
is harmonic, then $\varphi \circ f$ is harmonic only in
exceptional cases. Therefore one expects that harmonic mappings
from the disk into a strip domain behave quite differently from
harmonic mappings from the disk into a half-plane and that new
phenomena will be discovered in the study of harmonic maps. For
instance, it follows from  Lemma \ref{ClassSLemma} that
holomorphic functions  in plane do not increase hyperbolic
distances. In general, planar harmonic mappings do not enjoy this
property. On the other hand, we shall give here an additional
hypothesis under which the situation will change, in the plane as
well as in higher dimensions. It turns out that  the local uniform
boundedness  property, which we are going to define, has an
important role in our study.

For a domain $G \subset {\mathbb R}^n, n\ge 2, x, y \in G, $ let
$$r_{G}(x,y)= \frac{|x-y|}{\mathrm{min} \{ {d}(x),
  {d}(y) \}} \, \textrm{where} \, d(x)
= d(x, \partial G) \equiv \inf \{ |z-x|: z \in \partial G  \}\,.$$
If the domain $G$  is understood from the context,  we write  $r$
instead $r_G$. This quantity is used, for instance, in the study
of quasiconformal and quasiregular mappings, cf. \cite{vu2}. It is
a basic fact that \cite[Theorem 18.1]{v} for  $n\ge2, K\ge 1,
c_2>0$ there exists $c_1\in(0,1)$ such that whenever $f: G \to fG
$ is a 
$K$-quasiconformal mapping with $G, fG \subset {
\mathbb R}^n$ then $x,y \in G$ and $r_G(x,y) \le c_1 $ imply
$r_{fG}(f(x), f(y)) \le c_2.$ We call this property the local
uniform boundedness of $f$ with respect to $r_G\,.$ Note that
quasiconformal mappings satisfy the local uniform boundedness
property and so do quasiregular mappings under appropriate
conditions;   it is known that one to one mappings satisfying the
local uniform boundedness property    may  not be quasiconformal.
We also consider a weaker form of this property and say that $f: G
\to fG$ with $G, fG \subset { \mathbb R}^n$ satisfies the weak
uniform boundedness property
on $G$ (with respect to $r_G \,$) if there is a constant $c >0$ such
that  $r_G(x,y) \le 1/2$ implies $r_{fG}(f(x), f(y)) \le c\,.$
Univalent harmonic mappings fail to satisfy the weak  uniform
boundedness property as a rule, see Example \ref{myTh12} below.

We show  that if $f: G \to fG$ is  harmonic then   $f$ is
Lipschitz w.r.t. quasihyperbolic metrics on  $G$   and $fG$ if and
only if it satisfies the weak uniform boundedness  property; see 
Theorem \ref{m.res}. The proof is based on a higher
dimensional  version of the Schwarz lemma: harmonic maps  satisfy
the inequality  (\ref{harnack3}) below.  
An inspection of
the proof of 
Theorem \ref {m.res} shows that the class of harmonic
 mappings can be replaced by  $OC^1$ class defined by
(\ref{eq1}) (see Section 3 below)  and it leads  to
generalizations of the result;  see  
Theorem \ref{tmain}.


Another interesting  application is Theorem \ref{tqr}  which shows
that if $f$ is a harmonic $K$-quasiregular map such that the
boundary of the image is a continuum containing
at least two  points, then it is Lipschitz.
In  Subsection \ref{my24},  we study conditions under which a qc
mapping is quasi isometry with respect to the corresponding
quasihyperbolic metrics; see Theorems  \ref{thqh1} and
\ref{thqh2}. In particular, using a quasiconformal analogue of
Koebe's theorem, cf. \cite{ast.ge}, we give a simple proof  of the
following result, cf.  \cite{man,ma4}: if $D$ and $D'$ are proper
domains in $\mathbb{R}^2$ and $h : D \rightarrow D'$ is $K$-qc and
harmonic, then it is bi-Lipschitz with respect to quasihyperbolic
metrics on $D$ and $D'$.

The results in this paper may be generalized into various
directions. 
One direction is to consider weak 
continuous solutions of the $p$-Laplace equation
$$div(|\nabla u| ^{p-2}  \nabla u) = 0, \quad  1 < p < \infty, $$
so called $p$-harmonic functions.
Note that $2$-harmonic functions in the above sense  are
harmonic 
in the usual sense.

It seems that the case of the upper half space is  of
particular interest, cf. \cite{mk, ma4,ka,mak}.
In  Subsection
\ref{my25}, using Theorem 3.1 \cite{ksz} we prove  that if $h$ is a
quasiconformal $p$-harmonic mapping of the upper half space
$\mathbb{H}^n$ onto itself and $h(\infty)= \infty$, then $h$ is
quasi-isometry with respect to both the Euclidean and the
Poincar\'e distance.


\cc
\section{ Lipschitz property of harmonic maps w.r.t.
quasihyperbolic metrics}{}\label{sec2}

\subsection{Hyperbolic type metrics} \label{my21}
\bigskip

Let $B^n(x,r) = \{  z \in \mathbb{R}^n : |z-x|<r \}, S^{n-1}(x,r)=
\partial B^n(x,r)$ and let $\mathbb{B}^n, S^{n-1}$  stand for the
unit ball and the unit sphere in $ \mathbb{R}^n$, respectively.
Sometimes we write  $\mathbb{D}$ instead of  $\mathbb{B}^2$.
For a domain $G \subset \mathbb{R}^n$ let $\rho: G \to
(0,\infty)$ be a continuous function. We say that $\rho$ is a weight function
or a metric density if for every locally rectifiable curve
$\gamma$ in $G,$ the integral
$$ l_\rho(\gamma) = \int_{\gamma} \rho(x) ds$$
exists. In this case we call $ l_\rho(\gamma)$  the
$\rho$-length of $\gamma .$ A metric density defines a metric
$d_{\rho}: G \times G \to (0,\infty)$ as follows. For $a,b\in G,$
let
$$d_{\rho}(a,b)= \inf_{\gamma}  l_\rho(\gamma) $$
where the infimum is taken over all locally rectifiable curves in
$G$ joining $a$ and $b.$ For a fixed $a,b \in G\,,$
suppose that there exists a $d_\rho$-length minimizing curve
$\gamma:[0,1]\to G$ with $\gamma(0) =a, \gamma(1)=b$ such that
$$ d_{\rho}(a,b) = l_{\rho}(\gamma| [0,t]) + l_{\rho}(\gamma| [t,1]) $$
for all $t\in [0,1]\,.$ Then $\gamma$ is called a geodesic segment
joining $a$ and $b\,.$
It is an easy exercise to check that
$d_{\rho}$ satisfies the axioms of a metric. For instance, the
hyperbolic (or Poincar\'e) metric of the unit ball $\mathbb{B}^n$
and the upper half space $\mathbb{H}^n= \{x \in \mathbb{R}^n : x_n >0 \}$
are defined  in
terms of the densities $\rho(x)= 2/(1-|x|^2)$ and $\rho(x)=1/x_n\,,$
respectively. It is a classical fact that in both cases the length-minimizing
curves, geodesics, exist and that they
are circular arcs orthogonal to the boundary \cite{be}.
In both cases we have even explicit formulas for the distances:
\begin{equation} \label{myBn}
\sinh \frac{\rho_{\mathbb{B}^n}(x,y)}{2} =
\frac{|x-y|}{\sqrt{(1-|x|^2)(1-|y|^2)}}\,,
\quad x,y\in \mathbb{B}^n\,,
\end{equation}
and
\begin{equation} \label{myHn}
\cosh {\rho_{\mathbb{H}^n}(x,y)} =1+ \frac{|x-y|^2}{2 x_n y_n }\,,
\quad x,y \in \mathbb{H}^n\,.
\end{equation}

Because the hyperbolic metric is invariant under conformal
mappings, we may define the hyperbolic metric in any simply
connected plane domain by using the Riemann mapping theorem, see
for example  \cite{kl}. The Schwarz lemma may now be formulated by
stating that an analytic function from a simply connected domain
into another simply connected domain is a contraction mapping,
i.e. the hyperbolic distance between the images of two points is
at most the hyperbolic distance between
the points. The 
hyperbolic metric is often the natural metric in classical
function theory. For the modern mapping theory, which also
considers dimensions $n\ge 3\,,$ we do not have a Riemann mapping
theorem and therefore it is natural to look for counterparts of
the hyperbolic metric. So called hyperbolic type metrics have been
the subject of many recent papers. Perhaps the most important of
these metrics are the quasihyperbolic metric $ k_G$ and the distance
ratio metric $j_G$ of a domain $G \subset {\mathbb{R}^n}\,.$  They
are defined as follows.

\begin{nonsec}{The quasihyperbolic and distance ratio metrics}
\end{nonsec} 
Let $G \subset \mathbb{R}^n$ be a domain.
The quasihyperbolic metric $k_G$ is a particular case of the metric
$d_{\rho}$ when $\rho(x)=\frac{1}{d(x,\partial G)}$ (see
\cite{gp,gos,vu2}). It was proved in \cite{gos} that for given
$x,y\in G\,,$ there exists a geodesic segment of length $ k_G(x,y)$ joining them.
The distance ratio metric is defined for $x,y \in G$ by setting
 $$j_G(x,y)=\log(1+r_G(x,y))= \log (1+
 \frac{|x-y|}{\mathrm{min} \{ {d}(x),
  {d}(y) \}}) \,
$$
where $r_G$ is as in the Introduction. It is clear that
$$j_G(x,y) \leq r_G(x,y)\,.
$$
Some applications of these metrics are reviewed in \cite{vu3}.
The recent PhD theses \cite{li},
\cite{kl}, \cite{manthesis} study the quasihyperbolic geometry
or use it as a tool.

\begin{lem}\label{B} {\textrm (\cite{gp}, \cite[(3.4),
Lemma 3.7]{vu2})}
Let $G$ be a proper subdomain of  $ \mathbb{R}^n\,.$
\begin{enumerate}
\item[{\rm (a)}]   If $x,y \in G$ and $ |y-x| \leq d(x)/2$,
then   $k_G(x,y)\leq 2 j_G(x,y)\, .$
\item[{\rm (b)}] For $x,y \in G$ we have
$k_G (x,y) \geq j_G (x,y) \geq \log  \left (1+ \frac{|y-x|}{d(x)}\right )$.
\end{enumerate}
\end{lem}


\subsection{Quasiconformal and quasiregular maps} \label{myqcqr}

\begin{nonsec}{Maps of class ACL and ACL${}^n$}  \end{nonsec}  
For each integer $k=1,...,n$
we denote $R_k^{n-1}=\{ x\in R^n : x_k=0\}$. The orthogonal
projection $P_k : \mathbb{R}^n \to  \mathbb{R}_k^{n-1}$, is given by $P_k x=x- x_k\,e_k\,.$

Let  $I= \{ x\in \mathbb{R}^n : a_k\leq x_k \leq b_k\}$   be a
closed $n$-interval.
A  mapping  $f : I \rightarrow\mathbb{R}^m$
is said to be  absolutely continuous on lines (ACL) if $f$ is
continuous  and if $f$ is  absolutely continuous on almost
every  line  segment  in $I$, parallel  to the coordinate axes.
More precisely, if $E_k$  is the set of all  $x\in P_k I$  such
that the function $ t \mapsto u (x + t e_k)$ is not
absolutely continuous on $ [a_k,b_k]$,  then  $m_{n-1}(E_k)=0$ for
all $1\leq k\leq n$.

If $\Omega$ is an open set in $\mathbb{R}^n$,  a mapping  $f :
\Omega \rightarrow\mathbb{R}^m$  is absolutely continuous if $
f|I$  is ACL   for every  closed  interval  $ I \subset \Omega$.
If  $\Omega$ and  $\Omega'$  are domains in
$\mathbb{\overline{R}}^n$,   a homeomorphism $f: \Omega
\rightarrow \Omega' $ is called  ACL if  $f|\Omega\setminus
\{\infty,f^{-1}(\infty)\}$  is ACL.

If  $f : \Omega \rightarrow \mathbb{R}^m$ is ACL, then the partial
derivatives of  $f$  exist a.e. in $\Omega$,  and they are Borel
functions. We say that $f$ is ACL${}^n$ if the partials are locally integrable.

\begin{nonsec}{Quasiregular mappings}  \end{nonsec} 
Let $G\subset \mathbb{R}^n$ be a domain.
A mapping $f : G\to \mathbb{R}^n$ is said to
be {\it quasiregular\/} (qr) if $f$ is ACL${}^n$ and if there exists a constant
$K\ge 1$ such that
$$
  |f'(x)|^n\le KJ_f(x)\;,\;\;|f'(x)|=\max_{|h|=1}|f'(x)h|\;,
$$
a.e.\ in $G$. Here $f'(x)$ denotes the formal derivative of $f$ at $x\,,$
The smallest $K\ge 1$ for which this inequality is
true is called the {\it outer dilatation\/} of $f$ and denoted by $K\,_O(f)$.
If $f$ is quasiregular, then the smallest $K\ge 1$ for which the inequality
$$
 J_f(x)\le K\,l(f'(x))^n\;,\;\;l(f'(x))=\min_{ |h| =1}|f'(x)h|\;,
$$
holds a.e.\ in $G$ is called the {\it inner dilatation\/} of $f$ and denoted
by $K\, _I(f)$. The {\it maximal dilatation\/} of $f$ is the number
$K(f)=\max \{\,K\,_ I(f),\,K\,_ O(f)\,\}$. If $K(f)\le K$, then $f$ is said
to be $ K$--{\it
quasiregular\/} ($ K$--qr). If $f$ is not quasiregular, we set
$K\,_ O(f)=K\,_ I(f)=K(f)=\infty \,.$

Let  $ \Omega_1$ and $ \Omega_2$ be domains in $\mathbb{R}^n $
and fix $K\ge 1\,.$ We say
that a homeomorphism $ f : \Omega_1
\rightarrow \Omega_2$ is a $K$-quasiconformal (qc) mapping if it
is $K$-qr and injective.
Some of the standard references for qc and qr mappings are
\cite{ge2}, \cite{lv}, \cite{v}, and  \cite{vu2}.
These mappings generalize the classes of
conformal maps and analytic functions to Euclidean spaces.
The K\"uhnau handbook \cite{kuh} contains several reviews
dealing with qc maps.
It should be noted that various definitions for qc maps are studied
in \cite{v}. The above definition of $K$-quasiconformality is equivalent
to the definition based on moduli of curve families in \cite[p. 42]{v}.
It is well-known that qr maps are differentiable a.e., satisfy condition
(N) i.e. map sets of measure zero (w.r.t. Lebesgue's $n$-dimensional
measure) onto sets of measure zero. The inverse mapping of a $K$-qc mapping
is also $K$-qc. The composition of a $K_1$-qc and of a  $K_2$-qc map is
a $K_1 K_2$-qc map if it is defined.

\subsection{Examples} \label{my22}

We first show that, as a rule, univalent harmonic mappings  fail to satisfy
the local uniform boundedness property.

\begin{nonsec} {Example} \label{myTh12} \end{nonsec}
The univalent harmonic mapping $f : {\mathbb H}^2 \to  f({\mathbb
H}^2) \,,$ $f(z) =  \mathrm{arg}\, z + i \, \mathrm{Im} z,$ fails
to satisfy the local uniform boundedness property with respect to
$r_{{\mathbb H} ^2} \,.$

Let   $z_1= \rho e^{i \pi/4}$,     $z_2= \rho e^{i 3 \pi/4}$, $w_1
=f( z_1)$ and  $w_2 =f( z_2)$. Then $r_{{\mathbb H}^2}(z_1,z_2)=
2$  and   $r_{f{\mathbb H}^2}(w_1,w_2)=
\frac{\pi}{\sqrt{2}\,\rho}$  if  $\rho$ is small enough and  we
see that $f$  does  not satisfy the  local uniform boundedness
property.

\begin{nonsec} {Example}\label{1}\end{nonsec}
The univalent harmonic mapping $f: {\mathbb H}^2 \to  {\mathbb
H}^2\,,$ $f(z) = \mathrm{Re}\, z\,  \mathrm{Im}\, z + i \,
\mathrm{Im}  z,$ fails to satisfy the local uniform boundedness
property with respect to $r_{{\mathbb H} ^2} \,.$

For a harmonic mapping $f(z)= h(z)  + \overline{g(z)}\, ,$ we introduce
the following notation
$$
\lambda_f(z)= |h'(z)|-|g'(z)|\, ,  \quad \Lambda_f(z)=
|h'(z)|+|g'(z)| \,\,\,{\rm and} \quad \nu(z)= g'(z)/f'(z) .
$$


The following 
Proposition  shows   that a one to one
harmonic function satisfying the local uniform boundedness
property need not be quasiconformal. 
\begin{nonsec} {Proposition} \label{2} \end{nonsec}
The function $f(z)= {\mathrm{log}}(|z|^2) +2iy$ is a univalent
harmonic mapping and satisfies the local uniform boundedness
property, but $f$ is not quasiconformal on  $ V= \{z : x > 1,
0<y<1\}$.
\begin{pf}
It is clear that $f$ is harmonic in $\Pi^+ = \{ z :\, {\rm Re\,} z
> 0\}$. Next $f(z)= h(z) + \overline{g(z)}$, where $ h(z)=
{\mathrm{log} } z +z$ and $ g(z)= {\mathrm{log}} z - z$. Since
$h'(z)= 1 + 1/z$ and $g'(z)= -1 + 1/z$, we have $|\nu(z)| < 1$ for
$z\in \Pi^+$.


Moreover, $f$  is quasiconformal on every compact subset
$D\subset\Pi^+$ and  $\lambda_f$, $\Lambda_f$  are bounded from
above and below on $D$.  Therefore $f$ is a  quasi-isometry on $D$
and by Theorem \ref{m.res} below, $f$ satisfies the   local
uniform boundedness property  on $D$.

From now on we consider the restriction of $f$ to $ V= \{z =x + i
y : x > 1,  0<y<1\}$. Then   $fV= \{w =(u,v) : u > \log(1+ v^2/4),
0<v<2\} \,.$

We are going to show that:
\begin{enumerate}
\item[$\bullet$] $f$  satisfies the local uniform
boundedness   property,  but $f$ is not quasiconformal on $V$ .
\end{enumerate}

We see that $f$ is not quasiconformal on $V\,,$
because $|\nu(z)|\rightarrow 1$ as $ z\rightarrow \infty, z \in V\, .$

For  $s >1$,  define   $V_s= \{z : 1< x <s,  0<y<1\}$. Note      
that  $f$ is qc on  $V_s$   and  therefore  $f$  satisfies
the property of local uniform boundedness on $V_s$  for every $s>1$.      

We consider separately two cases.

Case A. $z\in V_4\,.$                           
If  $r>1$ is big enough,  then              
$d(z, \partial V_{r})=d(z,\partial V)$ and  
$d\big(f(z),\partial f(V_{r})\big)=d(f(z),\partial f(V))$
and therefore  $f$  satisfies the property of local uniform
boundedness on $V_4$  with respect to $r_V$.

Case B. It remains to prove that $f$ satisfies the property of local
uniform boundedness on $V\setminus V_4$  with respect to $r_V$.

Observe first that for $z, z_1 \in V$ and $|z_1|\geq
|z|\ge 1\,,$ we have the estimate
$$\log  \left ( \frac{|z_1|}{|z|}\right )\leq \frac{|z_1|}{|z|} -1\leq |z_1 -
z|,
$$
and therefore  for $z, z_1 \in V$              
\begin{equation}\label{ex1}
|f(z_1) - f(z)|  \leq 4 |z_1 - z| .
\end{equation}
 We write
$$  \partial V= [1,1+i] \cup A \cup B\,;   A= \{(x,0): x \ge 1\}\, ,
B= \{(x,1): x \ge 1\} \,.$$ Then
$$\partial (fV)= f(\partial V) \subset f[1,1+i] \cup (fA) \cup (fB)$$
and by the definition of $f$ we see that
$$fA=  \{(x,0): x \ge 0\}\,, \quad
fB = \{(x,2): x \ge \log 2\}, \quad f[1,1+i]
\subset [0, \log 2] \times [0,  2] \,.$$ Clearly for $w\in fV$
$$d(w) = \min \{ d(w, fA),d(w, fB),d(w, f[1,1+i]) \}\,, $$
and for ${\rm Re} w>1+ \log{2}$  and  $w \in fV\,,$  we find
\begin{equation} \label{eq1b}
  d(w) = \min \{ d(w, fA),d(w, fB) \}\,.
\end{equation}
For $z \in V  \setminus V_4$ we have ${\rm Re} f(z)\ge \log (16)
>1+  \log{2}$ and therefore, in view of the definition of $f$,
(\ref{eq1b})
 yields  $d(f(z))= 2 d(z)$. This together with  (\ref{ex1})   shows that $f$
satisfies the property of local uniform boundedness on $V\setminus
V_4$.
\end{pf}

\subsection{Higher dimensional version of Schwarz lemma} \label{my23}

Before giving  a proof of the higher dimensional version of the
Schwarz lemma for harmonic maps we first establish some notation.

Suppose that $h: \overline{B}^n(a,r)
\rightarrow\mathbb{R}^n$ is a continuous vector-valued function,
 harmonic on $ B^n(a,r)$, and let
$$M_a^*= \sup\{ |h(y)-h(a)| : y\in S^{n-1}(a,r) \}.
$$

Let $h=(h^1,h^2,\dots,h^n)$.  A modification of the estimate in
\cite[Equation (2.31)]{gil.trud} gives
$$r |\nabla h^k (a)| \leq n M_a^*\,,  \quad k=1,\dots,n.$$

We next extend this result to the case of vector valued functions.
See also  \cite{bu} and \cite[Theorem 6.16]{abr}.

\begin{lem}\label{le.sch}
Suppose that $h: {\overline{B}^n(a,r)}
\rightarrow\mathbb{R}^n$ is a continuous mapping, harmonic in
${B}^n(a,r)$.  Then
\begin{equation}\label{harnack2a}
r |h'(a)| \leq n M_a^*\,.
\end{equation}
\end{lem}

\begin{pf}
Without loss of generality, we may suppose that $a=0$ and
$h(0)=0$. Let
$$  K(x,y)=K_y(x)=  \frac{r^2-|x |^2}{n \omega_n \,r|x-y|^n},
$$
where $\omega_n$  is the volume of the unit ball $\mathbb{B}^n$ in ${\mathbb R}^n$.

Then
$$ h(x)=\int_{S^{n-1}(0,r)} K(x,t) h(t) d\sigma ,  \quad x\in  {B}^n(0,r),
$$
where $d\sigma$ is the $(n-1)$-dimensional  surface measure on $
S^{n-1}(0,r)$.                                                     

 A simple calculation yields

$$\frac{\partial }{ \partial x_j} K(x, \xi) =
\frac{1}{n \omega_n \,r} \left( \frac {-2x_j }{|x-\xi |^n}
 - n(r^2 - |x|^2) \frac{x_j - \xi_j}
{ |x-\xi |^{n+2}}\right).
$$

Hence, for $1 \leq j \leq n$, we have

$$\frac{\partial }{ \partial x_j} K(0, \xi) =
\frac{ \xi_j}{\omega_n\, r^{n+1}}\,.$$

Let $\eta \in S^{n-1} $ 
be a unit vector and $|\xi|=r$.
For  given  $\xi$, it is convenient  to write $K_\xi (x)=K(x,\xi)$
and consider   $K_\xi$ as a function of $x$.

Then
$$ K_\xi' (0) \eta= \frac{1}{\omega_n\, r^{n+1}}(\xi,\eta)\,.
$$

Since $|(\xi,\eta)|\leq |\xi| |\eta|= r$, we see that
$$|K_\xi' (0) \eta|\leq \frac{1}{\omega_n\, r^{n}} \,,
\textrm{
and therefore\,}\quad
|\nabla K^\xi (0)|\leq \frac{1}{\omega_n r^n}\,.
$$

This last inequality yields
$$|h'(0)({\eta})|\leq \int_{{S^{n-1}(a,r)} } |\nabla K^y(0)|\, |h(y)|\,
d\sigma(y)
\leq \frac{M_0^*\, n \omega_n r^{n-1} }{\omega_n r^n} =\frac{M_0^*
n}{r}
$$
and the proof is complete.
\end {pf}

Let $G \subset \mathbb{R}^n,$ be  a domain,
let $h: G \to {\mathbb R}^n$ be continuous. For  $x\in G\,$
let  $ B_x =B^n(x, \frac{1}{4}d(x))$  and
\begin{equation}\label{Mxdef}
M_x=\omega_h(x)=\sup\{ |h(y)-h(x)| : y\in B_x\}.
\end{equation}

If $h$  is a harmonic mapping, then
the inequality (\ref{harnack2a}) yields
\begin{equation}\label{harnack3}
\frac{1}{4} d(x) |h'(x)| \leq n \,\omega_h(x),\quad x\in G\,.
\end{equation}

We also refer to (\ref{harnack3}) as the inner gradient estimate.

\subsection{Harmonic quasiconformal  quasi-isometries} \label{my24}

For our purpose it is convenient to have the following lemma.

\begin{lem}\label{l.metr.1}
Let  $G$  and $G'$ be two domains in $\mathbb{R}^n$, and let $\sigma
$ and $\rho$  be two continuous metric densities on $G$  and $G'$,
respectively,
which define  the elements of length $ds= \sigma(z)|dz|$ and $ds=
\rho (w)|dw|$, respectively; and suppose that $f:G \rightarrow
G'$, \, is a $C^1$-mapping.\\
a)   If  there is a positive constant $c_1$  such that
$\rho(f(z))\, |f'(z)| \leq c_1 \, \sigma(z),$ $ z \in G\,,$ then $
d_{\rho} (f(z_2),f(z_1)) \leq c_1 \, d_{\sigma}(z_2,
z_1)$, \, $z_1, z_2\in G$.\\
b)  If  $f(G)=G'$  and  there is a positive constant $c_2$  such
that $\rho(f(z))\, l(f'(z)) \ge c_2 \, \sigma(z),$ $ z \in G\,,$
then $ d_{\rho} (f(z_2),f(z_1)) \ge c_2 \, d_{\sigma}(z_2, z_1)$,
\, $z_1, z_2\in G$
\end{lem}

The proof of this result is straightforward and it is left to the
reader as an exercise.
\begin{nonsec}{Pseudo-isometry
and a quasi-isometry}  \end{nonsec}  
Let $f$  be a map from a metric space $(M, d_M)$ into another
metric space $(N, d_N)$.
\begin{itemize}
\item   We say that $f$ is  a pseudo-isometry if there exist two
positive constants $a$ and $b$ such that for all $x, y\in M$,
$$
a^{-1} d_M(x,y) -b \leq d_N (f(x),f(y))\leq a d_M(x,y).
$$

\item   We say that $f$ is  a quasi-isometry or a bi-Lipschitz mapping
if there exists a positive constant $a \ge 1$ such that for all $x, y\in M$,
$$
a^{-1} d_M(x,y)  \leq d_N (f(x),f(y))\leq a d_M(x,y).
$$
\end{itemize}

For the convenience of the reader we begin our discusssion for
the unit disk case.

\begin{theo}{\label{n2}} Suppose that $h: {\mathbb D} \to {\mathbb R}^2$
is harmonic and satisfies the weak uniform boundedness property.
\begin{enumerate}
\item[{\rm (c)}] Then  $h :({\mathbb D}, k_{\mathbb D})
\to ( h({\mathbb D}), k_{h({\mathbb D}) })$ is Lipschitz.

\item[{\rm (d)}]  If, in addition, $h$ is a qc mapping,
then
$h :({\mathbb D}, k_{\mathbb D})
\to ( h({\mathbb D}), k_{h({\mathbb D}) })$
 is a quasi-isometry.
\end{enumerate}
\end{theo}
\begin{pf}
The part (d) is proved in \cite{ma4}.

For the proof of part (c) fix  $x\in \mathbb{D}$ and $y \in
B_x=B(x,\frac{1}{4} d(x) )$. Then $d(y) \geq \frac{3}{4}d(x)$ and
therefore $r(x,y)< 1/2$. By the hypotheses $ |h(y)-h(x)|\leq c\, d
(h(x))$.
The Schwarz lemma, applied to $B_x\,,$ yields in view of (\ref{Mxdef})
$$\frac{1}{4} d(x) |h'(x)| \leq 2 M_x  \leq 2 c \, d(h(x))
$$
The proof of part (c) follows from Lemma \ref{l.metr.1}.
\end{pf}

A similar proof applies 
for higher dimensions; the following result is a generalization of the part
$(c)$  of Theorem \ref{n2} .

\begin{theo}{}\label{m.res}
 Suppose that $G$ is  a proper subdomain of  $ \mathbb{R}^n$ and
$h : G \to \mathbb{R}^n$ is a harmonic mapping. Then the following
conditions are equivalent
\begin{enumerate}
\item[{\rm (1)}] $h$ satisfies the weak   uniform boundedness  property.
\item[{\rm (2)}]
$h :({G}, k_{G})
\to ( h({G}), k_{h({G}) })$
 is Lipschitz. 
\end{enumerate}
\end{theo}
\begin{pf}
Let us prove that $(1)$ implies $(2)$.

By the hypothesis  $(1)$  $f$ satisfies the weak   uniform
boundedness  property:    for every $x\in G$ and $t\in B_x$
\begin{equation}
|f(t) -f(x)|\leq c_2\, d(f(x))\,.
\end{equation}
This inequality together with  Lemma \ref{le.sch}
gives $d(x) |f'(x)|\leq c_3\, d(f(x))$  for every  $x\in G$. Now
an application of Lemma \ref{l.metr.1} shows that  $(1)$ implies
$(2)$.


It remains  to prove that  $(2)$ implies $(1)$.


 Suppose that $f$ is Lipschitz with the multiplicative constant
$c_2$. Fix $x,y \in G$ with $r_G(x,y)\leq 1/2$.  Then $|y-x| \leq
d(x)/2$ and therefore by Lemma \ref{B}
$$k_G(x,y)\leq 2 j_G(x,y)\leq 2\,r_G(x,y)\leq 1 .
$$
Hence $k_{G'}(fx,fy)\leq c_2$. Since  $j_{G'}(fx,fy)\leq
k_{G'}(fx,fy)\leq c_2 $,
we find  $j_{G'} (fx,fy) = \log ( 1 + r_{G'}(fx,fy)) \leq c_2$ and
therefore $r_{G'}(fx,fy)\leq e^{c_2} -1$.
\end{pf}

Since $f^{-1}$ is qc,  an application of  \cite[Theorem 3]{gos} to
$f^{-1}$ and  Theorem \ref{m.res}  give the following corollary:
\begin{cor}{}
Suppose that $G$ is a proper subdomain of  $ \mathbb{R}^n$,
$h:G\rightarrow hG$ is harmonic  and $K$-qc. Then $h :({G}, k_{G})
\to ( h({G}), k_{h({G}) })$ is a pseudo-isometry.
\end{cor}
In \cite[Example 11.4]{vu2} (see also \cite[Example 3.10]{vu1}),
it is  shown  that the analytic function  $f: {\mathbb D} \to G,
G= {\mathbb D} \setminus \{0 \}, f(z) = \exp((z+1)/(z-1))\,,$
$f({\mathbb D})=G,$ fails to be uniformly continuous as a map
$$f: ( {\mathbb D}, k_{ {\mathbb D}}) \to ( G, k_{ G})\,.$$
 Therefore bounded analytic functions do not  satisfy the weak
uniform boundedness property in general.   The situation will be
different for instance if the boundary of the image domain is a
continuum containing at least two points. Note that if $k_G$ is
replaced by the hyperbolic metric $\lambda_G$ of $G,$ then $f: (
{\mathbb D}, k_{ {\mathbb D}}) \to ( G, \lambda_{ G})\,$ is
Lipschitz.
\begin{theo}{}\label{tqr}
Suppose that $G \subset \mathbb{R}^n$,  $f :G\rightarrow
\mathbb{R}^n$ is  $K$-qr    and  $G'=f(G)$.  Let
$\partial G'$ be a continuum containing at least two distinct
points.
If $f$  is a harmonic mapping,  then
$f: ({G}, k_{G}) \to ({G'}, k_{G'})$
 is Lipschitz.
\end{theo}

\begin{pf}
Fix  $x\in G$ and let  $B_x=B^n(x,d(x)/4)$.
If $|y-x|\leq d(x)/4$, then $d(y)\geq 3d(x)/4$ and  therefore,
$$r_G(y,x)\leq\frac{4}{3} \frac{|y-x|}{d(x)}.
$$

Because  $j_G(x,y)=\log  (1+r_G(x,y))\leq r_G(x,y)$,  using Lemma
\ref{B}(a), we find
$$ k_G(y,x)\leq 2\,j_G(y,x)\leq 2/3 <1.
$$

By \cite[Theorem 12.21]{vu2} there exists a constant $c_2>0$ depending
only on $n$ and $ K$ such that
$$ k_{G'}(fy,fx)\leq c_2 \max \{k_G(y,x)^{\alpha}, k_G(y,x)\},
\alpha=K^{1/(1-n)}, $$
and hence,  using   Lemma \ref{B}(b) and $k_G(y,x)\leq 1$, we see that
\begin{equation}\label{qrB}
|fy-fx|\leq e^{c_2} d(fx),\quad \mbox{i.e.}  \quad M_x= \omega_f(x) \leq
e^{c_2}\, d(fx).
\end{equation}
By (\ref{harnack3}) applied to $B_x=B^n(x,d(x)/4)$, we have
$$\frac{1}{4} d(x) |f'(x)| \leq 2 M_x
$$
and therefore using the inequality  (\ref{qrB}),  we have
$$\frac{1}{4} d(x) |f'(x)| \leq 2\, c \,d(f(x)),$$
where $c=e^{c_2}$; and the proof follows from Lemma
\ref{l.metr.1}.                                               
\end{pf}

The first author has asked the following   Question (cf.
\cite{ma4}: Suppose that $G \subset \mathbb{R}^n$ is a proper
subdomain,  $f :G\rightarrow \mathbb{R}^n$ is harmonic $K$-qc and
$G'=f(G)$. Determine whether $f$ is a quasi-isometry
w.r.t. quasihyperbolic metrics on $G$ and $G'$.                          
This is true  for  $n=2$ (see Theorem  \ref{qi.planar} below). It
seems that one can modify the proof of Proposition 4.6 in
\cite{tw}  and show that this is true for the unit ball if $n \ge
3$ and $K<2^{n-1}$,
cf. also \cite{kam3}.
\subsection{Quasi-isometry  in planar case} \label{my25}

Astala and Gehring \cite{ast.ge} proved  a quasiconformal analogue
of Koebe's theorem, stated here as Theorem \ref{th.AGe}. These
concern the quantity

$$a_{f} (x) = a_{f,G} (x):= {\rm exp} \left(\frac{1 }{n|B_x|} \int_{B_x}
{\textrm{log}} J_f (z) dz \right),   x\in G, $$
associated with a quasiconformal mapping $f : G \rightarrow f(G)
\subset \mathbb{R}^n$; here $J_f$ is the Jacobian of f ; while
$B_x$ stands for the ball $B(x; d(x,\partial G)$; and $|B_x|$ for
its volume.
\begin{theo}{\cite{ast.ge}} \label{th.AGe} 
Suppose that G and $G'$  are domains in $R^n$: If  $f :
G\rightarrow G'$ is $K$- quasiconformal, then

$$\frac{1}{c}\frac{d (fx,\partial G')}{d(x,\partial G)}
\leq a_{f,G} (x) \leq c  \frac{d (fx,\partial G')}{d(x,\partial G)},
\quad   x\in G,$$
where  $c$ is a constant which depends only on $K$ and $n$.
\end{theo}

Let  $\Omega \in \mathbb{R}^n$  and  $\mathbb{R}^+ =[0,\infty)$.
If $f,g : \Omega \rightarrow \mathbb{R}^+ $  and there is a
positive constant $c$ such that

$$ \frac{1}{c}\, g(x)\leq f(x)\leq c\, g(x) \,, \quad x\in \Omega\,,$$
 we write
$f\approx g $  on $\Omega$.

Our next result  concerns  the quantity

$$E_{f,G} (x):= \frac{1 }{|B_x|} \int_{B_x}  J_f (z) dz\,,   x\in G, $$
associated with a quasiconformal mapping $f : G \rightarrow f(G)
\subset \mathbb{R}^n$; here $J_f$ is the Jacobian of $f $ ; while $B_x$ stands
for the ball $B(x, d(x,\partial G)/2$ and $|B_x|$ for its volume.

Define

$$A_{f,G}  =\sqrt[n]{E_{f,G}}\,.$$

\begin{theo}{}\label{thqh1}
Suppose  $f: \Omega\rightarrow \Omega'$ is a $C^1$ qc homeomorphism.
The following conditions are equivalent:

a)$f$  is bi-Lipschitz with respect to quasihyperbolic metrics on
$\Omega$ and $\Omega'  \,,$

b)  $\sqrt[n]{J_f} \approx  d_*/d\,,$

c)  $ \sqrt[n]{J_f}  \approx  a_f\,,$

d)  $\sqrt[n]{J_f}  \approx  A_f$,

\noindent
where  $ d (x)= d (x,\partial \Omega)$  and  $d_* (x)= d
(f(x),\partial \Omega')$.

\end{theo}
\begin{pf}
It is known that a) is equivalent to  b) (see for example
\cite{mm.unpub2}).

 In \cite{mm.unpub2}, using Gehring's result on  the distortion
property  of qc maps (see \cite{Ge}, p.383; \cite{v}, p.63),   the
first author  gives  short  proofs of  a new version of
quasiconformal analogue of Koebe's theorem;  it is proved that
$A_f \approx d_*/d$.

By Theorem \ref{th.AGe},  $a_f \approx  d_*/d$  and therefore  b)
is equivalent to  c). The rest of the proof is straightforward.
\end{pf}

If  $\Omega$ is planar domain  and $f$ a harmonic  qc map, then  we
proved that the condition d) holds.

The next theorem is a short proof of a recent result of V.
Manojlovic \cite{man},see also  \cite{ma4}.

\begin{theo} {}\label{qi.planar}
Suppose D and $D'$  are proper domains in $\mathbb{R}^2$. If  $h :
D \rightarrow D'$ is $K$- qc and harmonic, then it is bi-Lipschitz
with respect to quasihyperbolic metrics on D and $D'$.
\end{theo}

\begin{pf} Without loss of generality, we may suppose that $h$ is
preserving orientation. Let $z\in D$  and $h=f + \overline{g}$ be
a local representation of $h$ on $B_z$, where  $f$ and $g$  are
analytic functions on $B_z$, $\Lambda_h(z)= |f'(z)|+  |g'(z)|$,
$\lambda_h(z)= |f'(z)|- |g'(z)|$ and $k= \frac{K-1}{K +1}$.

Since  $h$  is  $K$- qc, we see that
\begin{equation}\label{est.jac}
(1-k^2) |f'|^2\leq J_h  \leq K |f'|^2
\end{equation}
 on $B_z$  and  since ${\textrm{log}} \, |f' (\zeta)|$ is harmonic,

$${\textrm{ log}} \, |f' (z)|= \frac{1 }{2|B_z|} \int_{B_z} \textrm{ log} \, |f' (\zeta)|^2 d\xi\, d\eta\,.$$

Hence, using  the right hand side of (\ref{est.jac}),  we find
\begin{eqnarray}
{\textrm{log}}\, a_{h,D} (z)\leq \frac{1}{2}\,
{\textrm{ log} }K + \frac{1 }{2|B_z|}
\int_{B_z} {\textrm{ log}} \, |f'(\zeta)|^2 d\xi\, d\eta\\
= {\textrm{ log}}\,\sqrt{K}\,|f'(z)|\,.
\end{eqnarray}

Hence,
$$a_{h,D} (z)\leq \sqrt{K}\,|f'(z)|$$
and in a similar way using  the left  hand side of (\ref{est.jac}),
we have

$$\sqrt{1-k^2}\, |f'(z)|\leq a_{h,D} (z)\,.$$
Now,  an application of
the Astala-Gehring result gives

$$\Lambda_h(z) \asymp \frac{d(hz,\partial D' )}{d(z,\partial D} \asymp \lambda_h(z)\,.$$

This pointwise result, combined with  Lemma \ref{l.metr.1}
(integration along curves), easily gives

$$k_{D'} (h(z_1), h(z_2)) \asymp k_D(z_1, z_2),\quad z_1, z_2 \in D\,.$$
\end{pf}

 Note that in \cite{man} the proof makes use of the interesting fact that
$ {\log} \frac{1}{J_h}$ is a
subharmonic function; but  we do not use it here.

Define $m_f(x,r)= \min \{ |f(x')- f(x)| : |x'-x|=r \}$.

Suppose that $G$ and $G'$  are domains in $\mathbb{R}^n$: If  $f :
G\rightarrow G'$ is $K$- quasiconformal; by the distortion property
we find  $m_f(x,r)\geq a(x) r^{1/\alpha}$. Hence, as in
\cite{kam3} and \cite{mm.unpub2} , we get:

\begin{lem}\label{Matelj}
If $f\in C^{1,1}$ is a $K-$ quasiconformal mapping defined in a
domain $\Omega\subset \Bbb R^n$ $(n\ge 3)$, then $$J_f(x)>0,\quad
x\in \Omega$$ provided that $K<2^{n-1}$. The constant $2^{n-1}$
is sharp.
\end{lem}

\begin{theo}{}\label{thqh2}
Under the hypothesis of the lemma, if $\overline{G}\subset
\Omega$, then 
$f$  is bi-Lipschitz with respect to Euclidean
and quasihyperbolic metrics on $G$  and $G'=f(G).$ 
\end{theo}

\begin{pf} Since $\overline{G}$ is compact $J_f$ attains minimum on $\overline{G}$ at
a point $x_0\in \overline{G}$. By Lemma \ref{Matelj},  $m_0=J_f
>0$ and therefore  since  $f\in C^{1,1}$ is a $K-$ quasiconformal, we conclude that
functions $|f_{x_k}|$, $1\leq k \leq n$  are bounded from above
and below on $\overline{G}$; hence  $f$ is bi-Lipschitz with
respect to Euclidean metric on  $G$.

By  Theorem  \ref{th.AGe}, we find   $a_{f,G} \approx d_*/d$,
where $ d (x)= d (x,\partial G)$ and $ d_* (x)= d (f(x),\partial
G')$. Since we have here $\sqrt[n]{J_f} \approx a_f$, we find
$\sqrt[n]{J_f} \approx d_*/d$ on $G$. An application  of Theorem
\ref{thqh1}
completes the proof. 
\end{pf}

\subsection{The upper half space  ${\mathbb H}^n.$} \label{my26}

Let  $\mathbb{H}^n$    denote the half-space in  $ \mathbb{R}^n$.
If  $D$ is a domain in  $ \mathbb{R}^n$,  by $QCH(D)$ we denote
the set of Euclidean harmonic quasiconformal mappings of $D$ onto
itself.

In particular if  $x\in  \mathbb{R}^3$,   we use notation
$x=(x_1,x_2,x_3)$ and we denote by $ \partial_{x_k} f= f'_{x_k}$
the partial derivative of $f$ with respect to $x_k$ .

A fundamental solution in space $ \mathbb{R}^3$ of the Laplace
equation is $ \frac{1}{|x|}$.
Let $U_0 = \frac{1}{|x+ e_3|}$, where $e_3=(0,0,1)$. Define
$h(x)=(x_1+ \varepsilon_1 U_0,x_2+ \varepsilon_2 U_0,x_3).$ It is
easy to verify that $h\in QCH (\mathbb{H}^3)$ for small values of
$\varepsilon_1$  and $\varepsilon_2$.

Using the Herglotz representation of   a nonnegative harmonic
function $u$  (see Theorem  7.24   and  Corollary  6.36
\cite{abr}), one can  get:

{\bf Lemma A}. If $u$ is a nonnegative harmonic function on a half
space $\mathbb{H}^n$, continuous up to the boundary with $u = 0$
on $\mathbb{H}^n$, then $u$ is (affine) linear.


 In  \cite{ma4},  the first author has outlined a proof
of the following result:

{\bf Theorem A}.   If  $h$ is  a quasiconformal  harmonic mapping
of the upper half space $\mathbb{H}^n$ onto itself and $h(\infty)=
\infty$, then $h$ is quasi-isometry with respect to both the
Euclidean and the Poincar\'e distance.

Note that  the outline of proof in \cite{ma4} can be justified by
Lemma A.

We show that  the analog statement of
this result holds  for $p$-harmonic vector functions (solutions of
$p$-Laplacian equations)  using the mentioned result  obtained in
the paper \cite{ksz}, stated here as:


{\bf  Theorem B}. If $u$ is a nonnegative $p$-harmonic function on a
half space $\mathbb{H}^n$, continuous up to the boundary with $u =
0$ on $\mathbb{H}^n$, then $u$ is (affine) linear.



\begin{theo} {}\label{th.som}
If   $h$ is  a quasiconformal  $p$-harmonic mapping of the upper
half space $\mathbb{H}^n$ onto itself  and $h(\infty)= \infty$,
then both $h: (\mathbb{H}^n, |\cdot|) \to  (\mathbb{H}^n,
|\cdot|)$
and $h: (\mathbb{H}^n, \rho_{\mathbb{H}^n}) \to
 (\mathbb{H}^n, \rho_{\mathbb{H}^n})$
are bi-Lipschitz where $\rho= \rho_{\mathbb{H}^n}$ is the
$Poincar\acute{e}$  distance.
\end{theo}
Since $2$-harmonic mapping are Euclidean harmonic  this result
includes Theorem A.


\begin{pf}
It suffices to deal with the case $n=3$ as the proof for the
general case is similar. Let $h=(h_1, h_2, h_3)$.


By Theorem B, we get  $h_3(x)=a x_3$, where $a$ is a positive
constant. Without loss of generality we may suppose that  $a=1$.

Since $h_3(x)=x_3$, we have $\partial_{x_3} h_3(x)=1$, and
therefore $|h'_{x_3}(x)|\geq 1$. In a similar way,
$|g'_{x_3}(x)|\geq 1$, where  $g=h^{-1}$. Hence,  there exists a
constant $c=c(K)$,
$$ |h'(x)|\leq c   \quad {\rm and} \quad    1/c  \leq  l(h'(x))\,.$$
Therefore partial   derivatives of $h$ and  $h^{-1}$ are bounded
from above;  and,  in particular, $h$ is  Euclidean bi-Lipschitz.

Since   $h_3(x)=x_3$,   $$\frac{|h'(x)|}{h_3(x)}\leq
\frac{c}{x_3};$$ and  hence, by Lemma \ref{l.metr.1},  $
\rho(h(a), h(b)) \leq c \rho (a, b)\,.$
\end{pf}


\bigskip


\section{Pseudo-isometry and  $OC^1(G)$}

In this  section,  we give  a sufficient condition  for  a qc
mapping  $f:G \rightarrow f(G)$  to be a pseudo-isometry w.r.t.
quasihyperbolic metrics on $G$ and $f(G)$.   First
we adopt the following notation.    

If  $V$ is a subset of $\mathbb{R}^n $ and $ u :V \rightarrow
\mathbb{R}^m$,  we define
$${\rm osc}_{V}u= \sup \{ |u(x) -u(y)| \, :  x, y \in V\}\, .$$

Suppose that $G \subset \mathbb{R}^n$  and   $B_x= B(x,d(x)/2)$. Let
$OC^1(G)$  
denote the class of $f \in C^1(G)$ such that
\begin{equation}\label{eq1}
d(x) |f'(x)|\leq c_1\,{\rm osc}_{B_x} f
\end{equation}
for every  $x\in G$. 
Similarly, let $SC^1(G)$ be the class of
functions $f \in C^1(G)$ such that
\begin{equation}
|f'(x)|\leq a r^{-1} \, \omega_f(x,r)  \quad  {\rm for\, all}\,\,
B^n(x,r)\subset G,
\end{equation}
where   $\omega_f(x,r)=\sup\{|f(y)-f(x)| : \,y \in B^n(x,r) \}$.\\

 The  proof of Theorem \ref{m.res}  gives the following more
general result:

\begin{theo}{} \label{tmain}
Suppose that $G \subset \mathbb{R}^n$,  $f :G\rightarrow G'$,
$f\in OC^1(G)$ and it satisfies the weak property of  uniform
boundedness  with a constant $c$ on $G$. Then
\begin{enumerate}
\item[{\rm (e)}]
$f: ({G}, k_{G}) \to ({G'}, k_{G'})$ is Lipschitz.
\item[{\rm (f)}] In addition, if  $f$ is   $K$-qc,
then $f$ is   pseudo-isometry
w.r.t. quasihyperbolic metrics on $G$ and  $f(G)$.
\end{enumerate}
\end{theo}
\begin{pf}
By the hypothesis  $f$ satisfies the weak property of  uniform
boundedness:   $ |f(t) -f(x)|\leq c_2\, d(f(x)$ for every $t\in
B_x$, that is
\begin{equation}
{\rm osc}_{B_x} f\leq c_2\, d(f(x))
\end{equation}
for every  $x\in G$. This inequality together with (\ref{eq1}) gives
$d(x) |f'(x)|\leq c_3\, d(f(x))$. Now an application of Lemma
\ref{l.metr.1} gives part  (e). Since $f^{-1}$ is qc,  an
application of  \cite[Theorem 3]{gos} on $f^{-1}$ gives part (f).
\end{pf}

 In order to apply the above method we
introduce  subclasses  of  $OC^{1}(G)$    (see, for example,
below (\ref{osc1})).                                                   


Let   $f:G\rightarrow G'$ be a $C^2$ function  and   $B_x=
B(x,d(x)/2).$
We denote by $OC^{2}(G)$ the class of        
functions which satisfy the following condition:
\begin{equation}\label{osc1}
\sup_{B_x} d^2(x) |\Delta f(x)|\leq c \,\, {\rm osc}_{B_x} f
\end{equation}
for every  $x\in G$.

If  $f\in OC^{2}(G)$,  then  by Theorem  3.9  in \cite{gil.trud},
applied to  $ \Omega=B_x$,

$$ \sup_{t \in B_x} d(t) |f'(t)|\leq C ( \sup_{t \in B_x} |f(t) -f(x)|
+ \sup_{t \in B_x} d^2(t) |\Delta f(t)|)$$

and hence  by  (\ref{osc1})
\begin{equation}
d(x) |f'(x)|\leq c_1\,{\rm osc}_{B_x} f
\end{equation}
for every  $x\in G$     and therefore $OC^{2}(G)\subset
OC^{1}(G)$.     

Now the following result follows from the previous theorem.

\begin{cor}{}
Suppose that $G \subset \mathbb{R}^n$ is a proper subdomain,
$f :G\rightarrow G'$ is  $K$-qc  and  $f$ satisfies  the
condition (\ref{osc1}).
Then $f: ({G}, k_{G}) \to ({G'}, k_{G'})$ is Lipschitz.
\end{cor}


We will now give some examples of classes of functions to which
Theorem \ref{tmain} is applicable.
  Let  $SC^2(G)$   denote the class of  $f \in
C^2(G)$ such that
$$ | \Delta f(x)|\leq a r^{-1} \sup\{|f'(y)| : \,y \in B^n(x,r)\},
$$
for  all  $ B^n(x,r)\subset G$,  where  $a$ is a positive constant.
Note that the class $SC^2(G)$   contains every function for
which $d(x) | \Delta f(x)|\leq a|f'(x)|$, $x\in G$. It is clear
that $SC^1(G) \subset OC^1(G)$ and by the mean value
theorem, $OC^{2}(G)\subset SC^2(G)$.
  For example,
in \cite{pa1} it is proved that  $SC^2(G)\subset SC^1(G)$ and that
the class $SC^2(G)$ contains  harmonic functions, eigenfunctions
of the ordinary Laplacian if $G$ is bounded, eigenfunctions of the
hyperbolic Laplacian if $G={\mathbb B}^n $ and thus our results
are applicable for instance to these classes.\\

\bigskip
{\bf Acknowledgements.} This work was initiated during the visit
of the first author to the Universities of Helsinki and Turku in
October 2005 and continued during the visit of the second author
to Belgrade in December 2006. The authors' research was supported
by the grant no. 8118525 the Academy of Finland. The first author
was partially supported by MNTRS, Serbia,  Grant No. 144 020. The
authors are indebted  
to the referee for valuable comments
and to O. Martio and S. Ponnusamy
for interesting discussions on this paper.

  \medskip
  \noindent
  {\bf  M. Mateljevi\'c }\\
  Faculty of Mathematics\\
  University of Belgrade\\
  Studentski trg 16\\
  11000   Beograd \\
  Serbia\\
  {\tt miodrag'at'matf.bg.ac.yu}\\

\medskip
\noindent
{\bf M. Vuorinen}\\
Department of Mathematics\\
FIN-20014 University of Turku \\
FINLAND\\
{\tt vuorinen'at'utu.fi}\\

\normalsize


\small

\begin{thebibliography}{1}

\bibitem {A} {\sc L. V. Ahlfors}:
{\em Conformal Invariants: Topics in Geometric Function Theory},
McGraw-Hill, New York, 1973.

\bibitem{akm} {\sc M. Arsenovi\'c, V. Koji\'c and M. Mateljevi\'c}:
{\it On Lipschitz continuity of harmonic
quasiregural mappings on the unit ball in $R^n$}, Ann. Acad. Sci.
Fenn. Math. Vol. {\bf 33}, (2008), 315--318. 

\bibitem{mak}   {\sc  M. Arsenovi\'c, V. Manojlovi\'c, and M.  Mateljevi\'c,}:
 {\it Lipschitz-type spaces and  harmonic mappings in the
space}, Ann. Acad. Sci. Fenn. 35 (2010), 1--9.

\bibitem{ast.ge} {\sc K. Astala and F. W. Gehring}:   { \it Quasiconformal analogues of
theorems of Koebe and Hardy-Littlewood},  Mich.Math.J. 32 (1985)
99-107.

\bibitem{abr}  {\sc S. Axler, P. Bourdon and W. Ramey}:
{\em Harmonic function theory}, Springer-Verlag, New York 1992.



\bibitem{be}
\textsc{A. F. Beardon:}
{\it The geometry of discrete groups}, Graduate Texts in Math. Vol 91,
Springer Verlag, Berlin -- Heidelberg -- New York, 1982.



\bibitem  {BH} {\sc D. Bshouty and W. Hengartner}:
\emph{Univalent harmonic mappings in the plane,}
  In: Handbook of Complex Analysis:
Geometric Function Theory, Vol. 2, (2005), 479--506,
      Edited by R. K\"uhnau (ISBN: 0-444-51547-X), Elsevier.

\bibitem  {bu} {\sc  B. Burgeth}: {\em A Schwarz lemma for harmonic
and  hyperbolic-harmonic functions in higher dimensions},
Manuscripta Math. {\bf 77} (1992), 283--291.

\bibitem{D} \textsc{P. Duren}: \emph{Harmonic mappings in the plane},
Cambridge University Press, 2004.



\smallskip
\bibitem{Ge} {\sc  F.W. Gehring}:
{\it Rings and quasiconformal mappings in  space},
Trans. Amer. Math. Soc. {\bf 103} , 1962, 353--393.



\bibitem{ge2} {\sc  F.W. Gehring}: {\it
 Quasiconformal mappings in Euclidean spaces.}
Handbook of complex analysis: geometric function theory. Vol. 2, 1--29, Ed. by R. K\"uhnau, Elsevier, Amsterdam, 2005.
\smallskip

\bibitem {gos} {\sc F.W. Gehring and B.G. Osgood}:
{\it Uniform domains and the quasi-hyperbolic metric}, J. Anal.
Math. {\bf 36}(1979), 50--74.

\bibitem  {gp} {\sc F. W. Gehring and B. P. Palka}:
{\em Quasiconformally homogeneous domains},
J. Anal. Math. {\bf 30} (1976), 172--199.

\bibitem {gil.trud} {\sc D. Gilbarg and   N. Trudinger}:
{\em Elliptic Partial
Differential  Equation of Second Order}, Second Edition, 1983.

\bibitem  {he} {\sc E. Heinz}:
{\em On one-to-one harmonic mappings}, Pacific J. Math. {\bf 9}(1959),
101--105.



\bibitem{ka} {\sc D. Kalaj}: {\it Quasiconformal and harmonic mappings between Jordan
domains}, Math. Z. {\bf 260:2}(2008),  237--252.


\bibitem{dkan}
{\sc D. Kalaj}: {\em Harmonic quasiconformal  mappings and Lipschitz
spaces},  Ann. Acad. Sci. Fenn. Math.  Math. {\bf 34:2} (2009), 475-485.






\bibitem {kam1} {\sc D. Kalaj and M. Mateljevi\'c}: {\em
Inner estimate and quasiconformal harmonic maps between smooth
domains}, J. Anal. Math. {\bf 100}, 117--132, (2006).

\bibitem{kam2} {\sc D.Kalaj and M. Mateljevi\'c}:  {\em Quasiconformal and harmonic
mappings between smooth Jordan domains},   Novi Sad J. Math. Vol.
{\bf 38} (2008), 147--156.

\bibitem{kam3}   {\sc  D. Kalaj and
M. Mateljevi\'c}:  {\em Harmonic quasiconformal self-mappings and
M\"obius transformations of the unit ball} , to appear in
Pacific J. Math.




\bibitem {kp}   {\sc D.
Kalaj and M. Pavlovi\'c}: {\em Boundary correspondence under quasiconformal harmonic diffeomorphisms of a half-plane}, Ann. Acad. Sci. Fenn. Math. {\bf 30} (2005), no. 1, 159--165.



\bibitem {kl} {\sc L. Keen and N. Lakic:}
{\it Hyperbolic geometry from a local viewpoint}.
London Mathematical Society Student Texts, 68. Cambridge University
Press, Cambridge, 2007.



\bibitem {ksz} {\sc T. Kilpel\"ainen, H. Shahgholian and X. Zhong}:
 {\em Growth estimates
through scaling for quasilinear partial differential equations},
Ann. Acad. Sci. Fenn. Math.  {\bf 32} (2007), 595--599.

\bibitem {kl} {\sc R. Kl\'en}: {\em On hyperbolic type metrics},
Ann. Acad. Sci. Fenn. Math. Diss.  {\bf 152} (2009), 1--49.



\bibitem {kuh} \textsc{R. K\"uhnau, ed.:}
{\it Handbook of complex analysis: geometric function theory}, Vol. 1-2.
  Elsevier Science B.V., Amsterdam, 2002, 2005.



\bibitem{lv}{\sc O. Lehto and K. I. Virtanen}:
{\em Quasiconformal Mappings in the Plane,} 2nd ed., Grundlehren Math.
Wiss., Band 126, Springer-Verlag, New York, 1973.



\bibitem {li} {\sc H. Lind\'en}: {\em Quasihyperbolic geodesics and
uniformity in elementary domains},
Ann. Acad. Sci. Fenn. Math. Diss.  {\bf 146} (2005), 1--50.


\comment{}
\bibitem{luz} {\sc A. Lyzzaik}:  {\em Local  properties of  Light
   Harmonic Mappings}, Canadian  J. Math.  {\bf 44}(1)(1992), 135--153.


\bibitem{manthesis} {\sc V. Manojlovi\'c}: {\em Moduli of Continuity of Quasiregular Mappings }, 2008, arXiv.math 0808.3241.


\bibitem{man} {\sc V. Manojlovic}: {\em  Bi-Lipschicity of quasiconformal harmonic
mappings in the plane}, Filomat,  {\bf 23} (2009), 85--89.


\bibitem{ma} {\sc O.Martio}: {\em
On harmonic quasiconformal mappings,} Ann. Acad. Sci. Fenn. Ser. A
I No. {\bf 425} (1968) 3-10.


\bibitem{ma3}
{\sc M. Mateljevi\'c}: {\em A version of Bloch  theorem for
quasiregular   harmonic mappings},  Rev. Roum.  Math. Pures Appl.
{\bf 47}(2002),  705--707.



\bibitem{ma4} {\sc M. Mateljevi\'c}:
{\em Distortion  of harmonic functions and harmonic quasiconformal
quasi-isometry},   Rev. Roum. Math. Pures Appl. {\bf 51}(2006),  711--722.


\bibitem{ma.koebe.ff} {\sc M. Mateljevi\'c}:
{\em Quasiconformal and quasiregular harmonic
analogues of Koebe's theorem and applications}, Ann. Acad. Sci.
Fenn.  {\bf 32}(2007), 301--315.

\bibitem{mm.unpub1} {\sc M. Mateljevi\'c}:
{\em On quasiconformal harmonic mappings}, unpublished manuscript,
2006.

\bibitem{mm.unpub2} {\sc  M. Mateljevi\'c}: {\em Lipschitz-type spaces,Quasiconformal and Quasiregular harmonic
mappings and Applications},  unpublished manuscript, 2008.


\bibitem{mk} {\sc M. Mateljevi\'c and M. Knezevi\'c}:
{\em  On the quasi-isometries   of harmonic  quasi-conformal
mappings},  J. Math. Anal. Appl. {\bf 334}(2007), 404--413.







\bibitem{ps} \textsc{D. Partyka and K. Sakan}: {\em
On bi-Lipschitz type inequalities for quasiconformal harmonic mappings},
 {Ann. Acad. Sci. Fenn. Math.} {\bf 32} (2007),  579--594.





\bibitem {pa1} {\sc M. Pavlovi\'c}:
{\em On  subharmonic behaviour of smooth functions},
Mat. Vesnik {\bf 48}(1996),  15--21.

\bibitem{pa2} {\sc M. Pavlovi\' c}:
{\em Boundary correspondence under harmonic quasiconformal
homeomorfisms of the unit disc}, Ann. Acad. Sci. Fenn., Vol 27,
2002, 365--372.

\bibitem {qiuvavu}  {\sc S. L. Qiu, M. K. Vamanamurthy, and M. Vuorinen}: {\em Some inequalities
for the Hersch-Pfluger distortion function}, J. Math. Anal. Appl.
Vol. {\bf 4} (1999), 2,  115-139
\comment{
\bibitem {riih} {\sc J. Riihentaus}: {\em A generalized mean value
inequality  for subharmonic functions},  Expo. Math. 19 (2001),
187--190.

\bibitem{s.yau} {\sc R. Schoen and S. T. Yau}:
{\em Lectures on Harmonic Maps\/}, Conf. Proc. and
Lect. Notes in Geometry and Topology, Vol. II, Inter. Press, (1997).
}
\bibitem{tw} {\sc L. Tam and  T. Wan}:
{\em On  quasiconformal   harmonic    maps},
Pacific J.  Math.  
 {\bf 182}(1998), 359--383.

\bibitem{v}    {\sc J.V\"{a}is\"{a}l\"{a}}:
{\em Lectures on n-Dimensional Quasiconformal Mappings},
Lecture Notes in Math. 229,  Springer-Verlag,  1971.


\bibitem{vu1} {\sc  M. Vuorinen}:
\emph{ Conformal invariants and quasiregular mappings.}
  J. Anal. Math. {\bf 45}(1985), 69--115.
\bibitem  {vu2} {\sc M. Vuorinen}: {\em Conformal Geometry
and Quasiregular Mappings}, Lecture Notes in Math. 1319,
Springer-Verlag, Berlin--New York, 1988.




\bibitem {vu3} {\sc M. Vuorinen}:
\emph{Metrics and quasiregular mappings.}
Proc. Int. Workshop on Quasiconformal Mappings and their Applications,
IIT Madras, Dec 27, 2005 - Jan 1, 2006, ed. by S. Ponnusamy, T. Sugawa and M. Vuorinen,
\emph{Quasiconformal Mappings and their Applications}, Narosa Publishing House, 291--325,
New Delhi, India, 2007.




\end{thebibliography}
\end{document}